\documentclass[12pt,a4paper]{article}
\usepackage{amsfonts}

\usepackage{graphicx}
\usepackage{amsmath}

\input{tcilatex}

\begin{document}

\title
{Rings of continuous functions, symmetric products, and
Frobenius algebras}

\author
{V.\,M.~Buchstaber [Bukhshtaber] and~E.\,G.~Rees\footnote{The present paper is based on the talk given by
E.~G.~Rees at the conference ``Kolmogorov and contemporary
mathematics,'' Moscow, 2003.}}
\date{19 March 2004}

\maketitle

\begin{abstract}

A constructive proof of the classical theorem of
Gel'fand and Kolmogorov (1939) characterising the image
of the evaluation map from a compact Hausdorff space $X$
into the linear space $C(X)^*$ dual to the ring $C(X)$
of the continuous functions on $X$ is given. Our approach 
to the proof enabled us to obtain a
more general result characterising the image of an
evaluation map from the symmetric products
$\operatorname{Sym}^n(X)$ into~$C(X)^*$. A similar
result holds if $X=\mathbb{C}^m$ and leads to explicit
equations for symmetric products of affine algebraic
varieties as algebraic sub-varieties in the linear space
dual to the polynomial ring. This leads to a better 
understanding  of the algebra of multi-symmetric
polynomials.

The proof of all these results is based on a formula
used by Frobenius in 1896 when defining higher
characters of finite groups. This formula had no further
applications for a long time; however, it 
occurred in several independent contexts during the last
fifteen years. The formula was used by A.~Wiles and R.~L.~Taylor
when studying representations and by H.-J.~Hoehnke and K.~W.~Johnson
and later by J.~McKay when studying finite groups. It
plays an important role in our work concerning
multi-valued groups. We describe several properties of
this remarkable formula. We also use it to prove a
theorem on the structure constants of Frobenius
algebras, which have recently attracted attention due to
constructions taken from topological field theory
and singularity theory. This theorem develops a result
of H.-J.~Hoehnke published in 1958. As a corollary, we
obtain a direct self-contained proof of the fact that
the 1-, 2-, and 3-characters of the regular
representation determine a finite group up to
isomorphism. This result was first published by H.-J.~Hoehnke 
and K.~W.~Johnson in 1992.

\end{abstract} 

{\large \bf Contents}

\bigskip
\begin{tabular}{llr} 
\S 1. & Introduction & 2 \\
\S 2. & Symmetric products & 6 \\ 
\S 3. & Properties of $n$-homomorphisms & 11\\
\S 4. & Frobenius algebras  & 15 \\
\S 5. & Appendix A. Proof of Mansfield's lemma & 20 \\
\S 6. & Appendix B. Algebra of multi-symmetric polynomials & 22 
\end{tabular}

\bigskip
\S 1. {\large \bf Introduction}
\bigskip

In~1939, A.~N.~Kolmogorov and I.~M.~Gel'fand published 
the paper `On rings of continuous functions on
topological spaces' (see~\cite{[11]}). The main result
of this paper identifies a compact Hausdorff space $X$
with the space of maximal ideals of the ring of
continuous functions on $X.$  Monographs and
textbooks containing basic  functional analysis usually
include this result. In modern terminology, it can be
stated as follows.

\medskip

{\bf Theorem~1} 
{\it Let $X$ be a compact Hausdorff
space. Then, for an appropriate topology on the space
of continuous complex-valued functions $C(X)$ on $X$,
the evaluation map 
$$
\mathcal{ E} : X\to\operatorname{Hom}(C(X),\mathbb{C}),
\qquad \mathcal{E}(x)\varphi=\varphi(x),
$$
is a homeomorphism onto the set of all ring
homomorphisms $C(X)\to\mathbb{C}$.}

This theorem is an analogue of Hilbert's Nullstellensatz: {\it If
$V$~is an affine algebraic variety with coordinate ring
$A=\mathbb{C}[x_1,\ldots ,x_n]/J$, where $J$~is a (radical) ideal
defining~$V$, then the evaluation map
$$
\mathcal{E} : V\to\operatorname{Hom}(A,\mathbb{C})
$$
defines an isomorphism between the variety~$V$ and the
set of all ring homomorphisms}.

Well-known reformulations of these theorems can be
obtained by noting that the set of all ring
homomorphisms $A\to\mathbb{C}$ can readily be identified
with the set (usually denoted by
$\operatorname{m-Spec}(A)$) of all maximal ideals
of~$A$.

We take the following point of view: the set of all ring
homomorphisms $f : A\to\mathbb{C}$ is an algebraic variety in
$\operatorname{Hom}(A,\mathbb{C})$ (regarded as a linear
space) defined by infinitely many equations 
$$
\{f(1)=1\text{ and }f(ab)=f(a)f(b)\text{ for every }a,b\in
A\} 
$$
involving the coordinate maps
$a : \operatorname{Hom}(A,\mathbb{C})\to\mathbb{C}$, $a(f)=f(a)$.
Since we deal with linear homomorphisms, it suffices to
consider the equations $f(1)=1$ and $f(a^2)=f(a)^2$,
where $a$ ranges over an additive basis in $A.$

The standard proof of the Kolmogorov--Gel'fand theorem
found in textbooks, argues by contradiction as follows. 

Let $I$ be a proper ideal in the function ring $C(X)$
such that there are no points $x\in X$ at which all
functions $\varphi\in I$ vanish. Then for any point
$x\in X$ there is a non-negative real-valued function
$\varphi_x\in I$ whose values are greater than~1 in some
neighbourhood $U_x$ of the point $x$. Using the fact
that the space $X$ is compact, we can  construct a
real-valued function $\varphi\in I$ whose values are
greater than~1 on the entire space~$X$. Such a function
is invertible, and hence a 
non-zero
constant function belongs to~$I$, that is, $I=C(X)$.

In \cite{[2]} we presented a constructive proof of our
theorem in~\cite{[3]} that characterises the symmetric products
$\operatorname{Sym}^n(X)$ as algebraic varieties
in~$C(X)^*$. In the case of $n=1$, this gives a proof of
the Gel'fand--Kolmogorov theorem. For a given ring
homomorphism $f : C(X)\to\mathbb{C}$, the proof constructs a
unique point $x\in X$ such that $f(\varphi)=\varphi(x)$
for any $\varphi\in C(X)$. In the literature
concerning this topic we found no proof of the
Gel'fand--Kolmogorov theorem similar to that presented
below.

Let us first consider the case in which $X$ is a finite
set. Then $C(X)\cong\mathbb{C}^n$, where $n$~is the number
of elements in~$X$. We choose a basis
$\{\delta_x:x\in X\}$ for~$C(X)$ such that
$\delta_x(x)=1$ and $\delta_x(y)=0$ if $x\ne y$. It is
clear that $1=\sum\delta_x$. Since
$\delta_x^2=\delta_x$, for a ring homomorphism $f$ we
have $f(\delta_x)^2=f(\delta_x)$, and hence
$f(\delta_x)=0$ or~$1$ for each $x$. On the other hand,
$1=f(1)=\sum_{x\in X}f(\delta_x)$, and hence there is a
unique point $x_0\in X$ such that $f(\delta_{x_0})=1$;
finally, using the expansion
$\varphi=\sum\varphi(x)\delta_x$, we see that
$f(\varphi)=\varphi(x_0)$ for any $\varphi\in C(X)$.

\smallskip Let us now show how to adapt this argument
to the case of a compact Hausdorff space~$X$.

\medskip
{\bf Definition 1}
{\it Let $K\subset X$ be a compact subset of a
topological space~$X$. A sequence of continuous functions
$\varphi_r : X\to[0,1]$ $(r\in\mathbb{N})$ is said to be {\it
enclosing\/} for $K$ if
\begin{itemize}
\item[1)]
$\varphi_r(\operatorname{Supp}(\varphi_{r+1}))=1$ for
any $r\in\mathbb{N}$,

\item[2)] $\varphi_r(x)=1$ for every $r\Leftrightarrow
x\in K$.
\end{itemize}}

It is clear that $\varphi_r\varphi_s=\varphi_s$ for
$s>r$.

\bigskip
\pagebreak
{\bf Example } 
Let $K$ be a compact subset of a
metric space $X$ with the metric
$d(\,\cdot\,{,}\,\cdot\,)$. Then the sequence of
functions
$$
\varphi_r(x)=\begin{cases} 0&\text{if}\;\;d(x,K)\ge1/r,
\\
1&\text{if}\;\;d(x,K)\le1/(r+1),
\\
(r+1)(1-rd(x,K))&\text{if}\;\;1/(r+1)\le d(x,x_0)\le1/r,
\end{cases}
$$
is enclosing for $K$.

The definition of an enclosing sequence for a compactum
$K$ obviously generalises to that of an {\it
enclosing net\/} (in other words, a {\it generalised
sequence\/}) of continuous functions
$\varphi_t : X\to[0,1]$ $(t\in T)$, where $T$~is some
directed set. The proof of the fact that an enclosing
net exists for any compactum~$K$ in a Hausdorff
space~$X$ is left to the reader. Moreover, to avoid
overloading our argument by technical details, we
speak in what follows about enclosing sequences, and 
use the fact that the assertions thus
obtained remain valid for enclosing nets.

\medskip

{\bf Lemma 2}
{\it Let $\{\varphi_r\}$ be an enclosing
sequence for a compactum~$K$. Then for any $r$ there is
an open neighbourhood $U_r$ of $K$ such that
$\varphi_r=1$ on~$U_r$.}

{\bf  Proof}
We set $U_r=\{x:\varphi_{r+1}(x)>0\}$. Then
$K\subset U_r$ and
$U_r\subset\operatorname{Supp}(\varphi_{r+1})$. Hence,
$\varphi_r=1$ on~$U_r$.

\medskip

{\bf Lemma 3} 
{\it Let $\{\varphi_r\}$ be an enclosing
sequence for a compactum $K$ and let $\psi : X\to\mathbb{R}$
be a function such that $\psi(x)=1$ for any $x\in U$,
where $U$~is an open neighbourhood of the compactum~$K$.
Then $(1-\psi)\varphi_r=0$ for~$r$ large enough.}

{\bf Proof} 
Let $x\notin U$. We choose an $r$ such that
$\varphi_r(x)=0$. Then the set $\varphi_r^{-1}[0,1)$ is
an open neighbourhood of~$x$ whose closure is disjoint
from~$U$. Let us cover the compact space $X\backslash U$
by finitely many neighbourhoods of this kind and denote
by $r_0$ the largest number among the numbers~$r$
corresponding to these neighbourhoods. Then
$(1-\psi)\varphi_{r_0}=0$, and hence
$(1-\psi)\varphi_s=(1-\varphi)\varphi_{r_0}\varphi_s=0$
for $s>r_0$.

\medskip

{\bf Lemma~4} 
{\it Let $f: C(X)\to\mathbb{C}$ be a ring
homomorphism and let $\varphi_r$~be a sequence of
functions in $C(X)$ such that
$\varphi_r\varphi_s=\varphi_s$ for any $r<s$. Then
either there is an $r_0$ such that $f(\varphi_r)=0$ for
$r\ge r_0$ or $f(\varphi_r)=1$ for any~$r$. }

{\bf Proof} We have $(f(\varphi_r)-1)f(\varphi_s)=0$
for any $r<s$. In this case, if $f(\varphi_r)\ne1$, then
$f(\varphi_s)=0$ for any $s>r$. This argument implies
the desired result.

\medskip

{\bf Definition~2} 
{\it Let $f : C(X)\to\mathbb{C}$ be a
ring homomorphism and let $\{\varphi_r\}$ be an
enclosing sequence for a compactum $K\subset X$. By the
{\it weight\/} of $K$ with respect to the sequence
$\{\varphi_r\}$ we mean the number
$w_f^\varphi(K)\in\{0,1\}$ equal to $f(\varphi_r)$ for
large values of~$r$.}

\medskip

{\bf Proposition 5} 
{\it Let $\{\varphi_r\}$ and
$\{\psi_r\}$ be two enclosing sequences for a given
compactum~$K$. Then $w_f^\varphi(K)=w_f^\psi(K)$.}

This defines the weight of a compactum, $w_f(K)$, which
does not depend on the choice of an enclosing sequence.

{\bf Proof} Suppose that $w_f^\varphi(K)=1$ and
$w_f^\psi(K)=0$. Then $f(\varphi_r)=1$ for any~$r$,
whereas $f(\psi_r)=0$ for any~$r$ greater than
some~$r_0$. Using Lemma~3 with $\psi=\psi_{r_0}$, we
find an $m>r_0$ such that
$\psi_{r_0}\varphi_m=\varphi_m$. Thus,
$f(1-\psi_{r_0})f(\varphi_m)=0$; however,
$f(1-\psi_{r_0})=1$, and hence $f(\varphi_r)=0$ for
$r>m$. The contradiction thus obtained proves the
Proposition. 

\medskip

{\bf Definition 3} 
{\it By the {\it support\/} of a
ring homomorphism $f : C(X)\to\mathbb{C}$ we mean the set
$S_f=\{x:w_f(x)=1\}$.}

\medskip

{\bf Proposition 6} 
{\it For any ring homomorphism $f$
the set $S_f$ consists of a single point.}

{\bf Proof}
We first assume that the set $S_f$ contains two distinct
points, say, $x$ and $y$. Let us choose an enclosing
sequence $\{\varphi_r\}$ for the set $\{x,y\}$. Then
$f(\varphi_r)$=1. For $r$ large enough we can choose
continuous functions $\psi_1$ and~$\psi_2$ such that
$\psi_1(x)=1$, $\psi_1(y)=0$, $\psi_2(x)=0$, and
$\psi_2(y)=1$,
$(\psi_1+\psi_2)^{-1}1\subset\varphi_r^{-1}1$, and
$\operatorname{Supp}\psi_1\cap\operatorname{Supp}\psi_2
=\emptyset$. Then $f(\psi_1)=f(\psi_2)=1$ by
construction, so, $f(\psi_1+\psi_2)=2,$ which is a
contradiction.

We now assume that the set $S_f$ is empty. Then for any point $x\in X$
there is a function $\varphi_x : X\to\mathbb{R}$ such that
$\varphi_x(x)=1$ and $f(\varphi_x)=0$. The open sets
$\{y:\varphi_x(y)>0\}$ cover the space~$X$.  Choosing a finite
sub-covering of this covering, we take the corresponding set of
functions $\varphi_1,\ldots ,\varphi_n$. Then the function
$\Phi=\varphi_1+\ldots +\varphi_n$ does not vanish at any point
of~$X$; however, $f(\Phi)=f(\varphi_1)+f(\varphi_2)+\ldots
+f(\varphi_n)=0$.  We see that, on the one hand,
$\displaystyle{f\left(\Phi\,\frac{1}{\Phi}\right)=f(1)=1}$, and, on the
other hand, $\displaystyle{f \left(\Phi \frac{1}{\Phi}\right)
=f(\Phi)f\left(\frac{1}{\Phi}\right)=0}$. This contradiction proves
Proposition~6.

Thus, to any ring homomorphism $f : C(X)\to\mathbb{C}$ we
assign a ring homomorphism $\widehat f : C(X)\to\mathbb{C}$,
$\widehat{f}(\varphi)=\varphi(x_0)$, where $S_f=\{x_0\}$.
To complete the proof of Theorem~1, it suffices to show
that $\widehat{f}=f$.

\medskip

{\bf Proposition 7} 
{\it Let $S_f=\{x_0\}$ and let
$\psi : X\to\mathbb{C}$ be a continuous function. Then
$f(\psi)=\psi(x_0)$.}

{\bf Proof}
Let $\{\varphi_r\}$ be an enclosing sequence for the
set $\{x_0\}$. We consider two cases.

We first assume that $\operatorname{Supp}\psi\cap
S_f=\emptyset$. In this case, there is an open
neighbourhood of the point $x_0$ on which the function
$\psi$ vanishes. Then $\psi\varphi_r=0$ for any~$r$
large enough, and at the same time we have
$f(\varphi_r)=1$. Hence, $f(\psi)=0$.

We now assume that $x_0\in\operatorname{Supp}\psi$. Let
us consider the function
$\theta_r=(\psi-\psi(x_0))\varphi_r$. We have
$|\theta_r|\le|\psi-\psi(x_0)|$, and $\theta_r$ vanishes
outside some neighbourhood of the point~$x_0$ because
of the factor~$\varphi_r$. Since the
function~$\psi$ is continuous, it follows that the
function $\theta_r$ tends to zero with respect to the
sup norm, and hence $f(\theta_r)\to0$ as $r\to\infty$.
On the other hand,
$f(\theta_r)=(f(\psi)-\psi(x_0)f(1))f(\varphi_r)\to
(f(\psi)-\psi(x_0))w_f(x_0)$ for~$r$ large enough.
Hence, $f(\psi)=\psi(x_0)$. This proves Proposition~7.

\bigskip

\S 2. {\large \bf Symmetric products }

\bigskip
We recall that by the symmetric product of a space~$X$
one means the quotient  space
$$
\operatorname{Sym}^n(X)=X^n/S_n=
\bigl\{(x_1,\ldots ,x_n):(x_{\sigma(1)},\ldots ,x_{\sigma(n)})
\sim(x_1,\ldots ,x_n),\ \sigma\in S_n\bigr\},
$$
where $S_n$ is the group of all permutations of a set
with $n$ elements.

The continuous functions on $\operatorname{Sym}^n(X)$
correspond  exactly to the continuous functions
$f : X^n\to \mathbb{C}$ invariant under all
permutations of the coordinates, that is, the 
symmetric functions.

Let us consider an analogue of the evaluation map 
$$
\mathcal{ E} :  \operatorname{Sym}^n(X)\to 
\operatorname{Hom}(C(X),\mathbb{C}), \qquad
\mathcal{E}(x_1,\ldots ,x_n)(\varphi)=\varphi(x_1)+\ldots +\varphi(x_n).
$$
We shall describe the image of this map by using
equations. These equations are given by formulas which
were first used by G.~Frobenius \cite{[9]}, \cite{[10]} and,
more recently, by a number of authors, including
A.~Wiles~\cite{[19]}, R.~L.~Taylor~\cite{[18]}, H.-J.~Hoehnke and
K.~W.~Johnson~\cite{[13]}, R.~Rouquier~\cite{[17]},
and L.~Nyssen~\cite{[15]}. The formulae play an important
role in the theory of multi-valued groups
(see~\cite{[4]}--\cite{[6]}). We follow the approach
developed in~\cite{[3]}.

Let $A$ be an associative algebra with unit over the
field $\mathbb{C}$ of complex numbers and let $f : A\to\mathbb{C}$ 
 be a linear tracial map (that is, $f(ab)=f(ba)$ for every $a,b \in
A$). We introduce linear maps $\Phi_n(f) : A^{\otimes
n}\to\mathbb{C}$ by setting $\Phi_1(f)=f$,
$\Phi_2(f)(a_1\otimes a_2)=f(a_1)f(a_2)-f(a_1a_2)$, and
further by recurrence,
$$
\Phi_{n+1}(f)(a_1\otimes a_2\otimes\ldots \otimes
a_{n+1}) =f(a_1)\Phi_n(f)(a_2\otimes\ldots \otimes
a_{n+1})$$
$$ -\Phi_n(f)(a_1a_2\otimes
a_3\otimes\ldots \otimes a_{n+1})
-\ldots -\Phi_n(f)(a_2\otimes a_3\otimes\ldots \otimes
a_1a_{n+1}).
$$
We note that a ring homomorphism $f : A\to\mathbb{C}$
satisfies the conditions $f(1)=1$ and
$\Phi_2(f)\equiv 0$.

\medskip
{\bf Definition 4} 
{\it By a {\bf Frobenius
$n$-homomorphism\/} we mean a linear homomorphism
$f : A\to\mathbb{C}$ satisfying the conditions $f(1)=n$ and
$\Phi_{n+1}(f)\equiv0$.}

Our choice of name for the above homomorphisms is
explained by the fact that the formula for the defining
recursion first 
arose in the papers of G.~Frobenius~\cite{[9]}
and~\cite{[10]} in the case of group algebras of finite
groups. For instance, the following result (in our notation) 
was obtained in~\cite{[9]}.

Let $G$ be a finite group and let $A=\mathbb{C}G$ be its
group algebra. Then the character $\chi : G\to\mathbb{C}$ of
any $n$-dimensional linear representation of the
group~$G$ can be extended to a linear homomorphism
$\chi : A\to\mathbb{C}$ such that $\chi(1)=n$ and
$\Phi_{n+1}(\chi)\equiv0$.

A generalization of the Gel'fand--Kolmogorov theorem to
the case of symmetric products is given by the following
result.

\medskip

{\bf Theorem~8} 
{\it Let $X$ be a compact Hausdorff
space. Then the image of the map 
$$
\mathcal{E} : \operatorname{Sym}^n(X)\to
\operatorname{Hom}(C(X),\mathbb{C}) 
$$
is exactly the subspace of all Frobenius
$n$-homomorphisms, that is, it is given by the equations
$f(1)=n$ and $\Phi_{n+1}(f)\equiv0$.}

A more detailed statement and a constructive proof of
this theorem which uses the above technique of enclosing
sequences (and enclosing nets in the case of general
compact Hausdorff space) can be found in~\cite{[2]}.

Similarly to Definition~4, one can introduce the notion
of Frobenius $n$-ho\-mo\-mor\-phisms $f : A\to B$, where
$B$~is an arbitrary commutative algebra. If $B$ has no
zero divisors, then the Frobenius $n$-homomorphisms have
properties which are important for our purposes.
Therefore, in what follows we assume that $B$~is an
integral domain. In~\cite{[3]}, the proof of Theorem~8 is
obtained as a consequence of a general result
characterizing  Frobenius $n$-homomorphisms $f : A\to B$.

Let $A$ be a commutative algebra. We denote by $S^nA$
the symmetric subalgebra of~$A^{\otimes n}$. Every
element $a\in S^nA$ can be represented in the form 
$$
{\bf a}=\sum_{\sigma\in S_n}a_{\sigma(1)}\otimes\ldots\otimes
a_{\sigma(n)} 
$$
and the product in $S^nA$ in the form
$$
{\bf ab}=\sum_{\sigma_1,\sigma_2\in S_n}
a_{\sigma_1(1)}b_{\sigma_2(1)}\otimes\ldots\otimes
a_{\sigma_1(n)}b_{\sigma_2(n)}.
$$

\medskip

{\bf Theorem~9} 
{\it A linear map $f:A\to B$ is a
Frobenius $n$-homomorphism if and only if $f(1)=n$ and
the restriction of the homomorphism
$\Phi_n(f) : A^{\otimes n}\to B$ to $S^nA$ gives a ring
homomorphism 
$$
\frac1{n!}\Phi_n(f) : S^nA\to B.
$$}

We now present a combinatorial result heavily used in
the proof of this theorem.

Let $X$ be a finite set and let $\mathcal{ P}(X)$ be the
free Abelian group generated by the set of all
partitions of the set~$X$. We recall that every
permutation~$\sigma$ of the set~$X$ defines a partition
of~$X$ given by the orbits of the action of the subgroup
generated by~$\sigma$. All the permutations defining the
same partition have the same sign; for a given partition
$\pi$ we denote by $\epsilon(\pi)$ the corresponding
sign and by $n(\pi)$~the number of permutations
generating~$\pi$. We set 
$$
\chi(X)=\sum_\pi\epsilon(\pi)n(\pi)\pi\in\mathcal{ P}(X),
$$
where the sum is taken over all partitions of the set
$X$.

Let $\pi_1$ and $\pi_2$ be partitions of sets $X$ and
$Y$, respectively. Then a natural partition $\pi_1\pi_2$
of the disjoint union $X\sqcup Y$ is defined. Thus, the
element $\chi(X)\chi(Y)\in\mathcal{ P}(X\sqcup Y)$ is well
defined.

If $g : X\to Y$ is a surjection and $\pi$ is a
partition of the set~$Y$, then one can take the
preimages of the parts of~$\pi$ and obtain a partition
$g^*\pi$; thus, an induced homomorphism $g^* : \mathcal{
P}(Y)\to\mathcal{ P}(X)$ is defined.

If the maps $i_1 : X\to Z$ and $i_2 : Y\to Z$ are
embeddings and if $Z=i_1(X)\cup i_2(Y)$, then we say
that $(Z,i_1,i_2)$ is an {\it amalgamated union\/} of
the sets $X$ and~$Y$. There is a map $q : X\sqcup Y\to
Z$, which is surjective.

The following rather unexpected purely combinatorial
result was obtained in~\cite{[3]} using 
polynomials which are solutions of a hypergeometric
differential equation.

\medskip

{\bf Lemma~10} 
{\it $$
\sum q^*\chi(Z)=\chi(X)\chi(Y),
$$
where the sum is taken over all distinct amalgamated
unions of $X$ and~$Y$, including the disjoint union.}

The main steps in the proof of Theorem~9 are as
follows.

1. Let $X=(a_1,\ldots ,a_n)$, where $a_k\in A$. Using the
homomorphism  $f : A\to B$, we define a homomorphism
$$
f : \mathcal{ P}(X)\to B, 
$$
whose value on the partition
$\pi=(P_1,\ldots ,P_k)\in\mathcal{ P}(X)$, where
$P_i=(a_{i_1},\ldots ,a_{i_q})$, is equal to 
$$
f(\pi)=\prod_{i=1}^kf(a_{i_1}\ldots a_{i_q}). 
$$
Then 
$$
f(\chi(X))=\Phi_n(f)(a_1,\ldots ,a_n).
$$
2. Let $X=(a_1,\ldots ,a_n)$ and $Y=(b_1,\ldots ,b_n)$ be
disjoint sets. It follows from the definitions that 
$$
f(\chi(X)\chi(Y))=f(\chi(X))f(\chi(Y)).
$$
3. Let $f : A\to B$ be an $n$-homomorphism. Using
Lemma~10 we obtain 
$$f(\chi(X)\chi(Y))=\sum_{\sigma \in S_n} f(\chi(X\sqcup_\sigma Y)). 
$$
Here $Z=X\sqcup_\sigma Y$ stands for an amalgamated
union such that $Z$ consists of $n$ elements; in this
case, the maps $i_1$ and~$i_2$ are one-to-one, and
therefore $Z$ is defined by some permutation~$\sigma$.

4. The last step of the proof is to use the relation 
$\Phi_{n+1}(f)\equiv 0$ to establish the
equality 
$$
\sum_{\sigma\in S_n}\Phi_n(f)(a_1b_{\sigma(1)},
\ldots ,a_nb_{\sigma(n)}) =n!\,\Phi_n(f)({\bf ab}).
$$
Let us now denote by $\Phi_n(A)$ the set of all
Frobenius $n$-homomorphisms of the ring~$A$ to the
field~$\mathbb{C}$ of complex numbers. By construction,
$\Phi_n(A)$~is an algebraic subvariety of the linear
space $A^*=\operatorname{Hom}(A,\mathbb{C})$ with the
coordinates $a : A^*\to\mathbb{C}$, where $a(f)=f(a)$. We set
$$
\Phi_n(\mathbb{C}[u_1,\ldots ,u_m])=\Phi_n(m).
$$

\medskip

{\bf Theorem 11}
{\it The map 
$$\mathcal{E} : \operatorname{Sym}^n(\mathbb{C}^m)\to
\operatorname{Hom}(\mathbb{C}[u_1,\ldots ,u_m],\mathbb{C}),$$
$$\mathcal{E}(x_1,\ldots ,x_n)(p)=p(x_1)+\ldots+p(x_n), $$
defines a homeomorphism 
$$ \mathcal{E} : \operatorname{Sym}^n(\mathbb{C}^m)\to\Phi_n(m).$$}
Using this result, in Appendix~B we present a
description of an embedding 
$$
\operatorname{Sym}^n(\mathbb{C}^m)\subset \mathbb{C}^N, \quad
N= \left(\begin{array}{c}n+m\\ n \end{array}\right) - 1 ,
$$
in the context of the theory of multi-symmetric
polynomials.

{\bf Proof} The proof of the fact that the evaluation map
$\mathcal{E}$ is an embedding is immediate (see~\cite{[3]}). Let
$f : \mathbb{C}[u_1,\ldots ,u_m]\to\mathbb{C}$ be a Frobenius
$n$-homomorphism. Then, by Theorem~9, the map
$$\frac{1}{n!}\Phi_n(f) : S^n(\mathbb{C}[u_1,\ldots ,u_m])\to\mathbb{C}$$
is a ring homomorphism. In what follows, the symmetric algebra
$S^n(\mathbb{C}[u_1, \ldots ,u_m])$ can be identified with the algebra
of polynomial functions on the algebraic variety
$\operatorname{Sym}^n(\mathbb{C}^m)$, that is, with the algebra of
multi-symmetric polynomials. Then, by Hilbert's Nullstellensatz, there
is a set of points
$(x_1,\ldots ,x_n)\in\operatorname{Sym}^n(\mathbb{C}^m)$ such that for
$a=\sum_{\sigma\in S_n}p_{\sigma(1)}\otimes\ldots \otimes p_{\sigma(n)}$
we have the formula
$$ a(x_1,\ldots ,x_n)=\sum_{\sigma\in S_n}
p_{\sigma(1)}(x_1)\otimes\ldots \otimes
p_{\sigma(n)}(x_n).
$$
Let us take for $a$ the element 
$$
a=p\otimes1\otimes\ldots \otimes1+1\otimes
p\otimes\ldots \otimes1
+\ldots +1\otimes1\otimes\ldots \otimes p,
$$
where $p\in\mathbb{C}[u_1,\ldots ,u_m]$. Then
$$
\frac1{n!}\Phi_n(f)(p)=\sum_{k=1}^np(x_k).
$$
On the other hand, using the fact that the homomorphism
$\Phi_n(f)$ is linear and symmetric and applying the
formula  
$$
\Phi_n(f)(p\otimes1\otimes\ldots \otimes1)
=f(p)(f(1)-1)\ldots(f(1)-(n-1))
$$
(see~\cite{[3]}), we obtain 
$$
\frac1{n!}\Phi_n(f)(a)=f(p)
$$
because $f(1)=n$, that is, $f(p)=\sum_{k=1}^np(x_k)$,
and hence $f=\mathcal{E}(x_1,\ldots ,x_n)$. This completes
the proof of Theorem~11.

\medskip

{\bf Theorem 12} 
{\it Let $A$ be a finitely generated
commutative algebra and let $V$~be the affine algebraic
variety $\operatorname{m-Spec}(A)$. Then the map 
$$
\mathcal{E} : \operatorname{Sym}^n(V)\to
\operatorname{Hom}(A,\mathbb{C}),
\quad
\mathcal{E}(f_1,\ldots ,f_n)(a)=f_1(a)+\ldots+f_n(a),
$$
where $f_k : A\to\mathbb{C}$, $k=1,\ldots ,n$, are ring
homomorphisms (whose kernels are the maximal ideals which 
are points of $V$), defines a homeomorphism 
$$
\mathcal{E} : \operatorname{Sym}^n(V)\to\Phi_n(A).$$}
For a detailed proof of this result, see \cite{[3]}.

The proof of Theorem~8 using Theorem~9 can be carried
out by the scheme of the proof of Theorem~11. One must
also use that for a compact Hausdorff
space $X$ the ring $C(X)^{\otimes n}$ is dense in the
ring $C(X^n)$ by the Stone--Weierstrass theorem (for the
details, see~\cite{[16]} Proposition~1.10.21), and
therefore the ring $S^n(C(X))$ is dense in the ring
$C(\operatorname{Sym}^n(X))$. This enables one to apply
the Gel'fand--Kolmogorov theorem and to assign a point
of the space $\operatorname{Sym}^n(X)$ to the ring
homomorphism $\frac1{n!}\Phi_n(f) : S^n(C(X))\to\mathbb{C}$.

\bigskip

\S 3. {\large \bf Properties of $n$-homomorphisms}

\bigskip
The results of the previous section show that 
Frobenius $n$-homomorphisms of commutative algebras have
important applications. Thus, the study of the
properties of these homomorphisms in this special case
arose naturally. It turns out that in this investigation
it is useful to apply another description of the
defining equations, which was given in~\cite{[4]},~\cite{[5]}.

Let $A$ be a commutative algebra. Then the
following assertions hold.

1. The value
$\Phi_n(a,\ldots ,a):=\Phi_n(f)(a\otimes\ldots \otimes a)$,
where $a\in A$, is equal to the determinant of the
matrix 
$$
\left(\begin{array}{cccccc}f(a) & 1 & 0 & 0 & \ldots&0
\\
f(a^{2}) & f(a) & 2 & 0 & \ldots&0
\\
f(a^{3}) & f(a^{2}) & f(a) & 3 & &0
\\
\vdots&\vdots&\ddots&\ddots&\ddots&
\\
\vdots&\vdots& & &f(a) & n-1
\\
f(a^{n}) & f(a^{n-1}) & \ldots&\ldots& f(a^{2}) & f(a)\end{array}\right).
$$

2. The value
$\Phi_n(f)(a_1,\ldots ,a_n)=\Phi_n(f)(a_1\otimes\ldots \otimes
a_n)$ is a multilinear function of the variables
$a_1,\ldots ,a_n$ and can therefore be obtained from
$\Phi_n(a,\ldots ,a)$ by the standard polarization
procedure.

3. Let us represent the permutation $\sigma\in S_n$ in
the form of a product of disjoint cycles of total length
$n$, say, $\sigma=\gamma_1\ldots\gamma_r$. Let
$\gamma=(i_1,\ldots ,i_m)$ be a cycle. We set
$f_\gamma(a_1,\ldots ,a_n)=f(a_{i_1}\ldots a_{i_m})$. Then
the following formula holds:
$$
\Phi_n(f)(a_1,\ldots ,a_n)=\sum_{\sigma\in
S_n}\varepsilon(\sigma) f_{\gamma_1}(a_1,\ldots ,a_n)\ldots
f_{\gamma_r}(a_1,\ldots ,a_n),
\hspace{1cm} (3.1)
$$
where $\varepsilon(\sigma)$ is the sign of the
permutation $\sigma$.

For instance, 
$
\Phi_3(f)(a_1,a_2,a_3)= $
$$f(a_1)f(a_2)f(a_3)-f(a_1)f(a_2a_3)
-f(a_2)f(a_1a_3)-f(a_3)f(a_1a_2)+2f(a_1a_2a_3).
$$

We note that formula (3.1) holds also in the case
of tracial homomorphisms~$f$ of a non-commutative
algebra $A$. This formula arose in~\cite{[7]} in
the case of matrix algebras~$A$ and the trace
homomorphism~$f$.

It is often simpler to work with the formula
(3.1) in its ``diagonal'' form and then to use
polarization to obtain general formulae. For
instance, 

\noindent
$\Phi_3(f)(a,a,a)=s_1^3-3s_1s_2+2s_3$, where
$s_1=f(a)$, $s_2=f(a^2)$, and $s_3=f(a^3)$.
We use the notation $s_k$ to stress the relationship
with the classical Newton formula expressing the
elementary symmetric functions $e_n=\sum t_{i_1}\ldots
t_{i_n}$ of the commuting variables $t_1,t_2,\ldots $ (the
sum is taken over all sets 

\noindent
$I_n=(i_1<\ldots <i_n)$) as
polynomials in the power sums $s_k=\sum t_i^k$, 
$$
n!\,e_n=\det \left(\begin{array}{cccccc} s_1 & 1 & 0 & 0 & \ldots&0
\\
s_2 & s_1 & 2 & 0 & \ldots&0
\\
s_3 & s_2 & s_1 & 3 & &0
\\
\vdots&\vdots&\ddots&\ddots&\ddots&
\\
\vdots&\vdots& & &s_1 & n-1
\\
s_n & s_{n-1} & \ldots&\ldots& s_2 & s_1
\end{array}\right).
$$

Let us now describe the development of this
relationship, denoting by $F_n$ the polynomial in the
variables $s_k$, $k=1,\ldots ,n$, which becomes equal to
$\Phi_n(f)(a,\ldots ,a)$ when making the substitution
$s_k=f(a^k)$. The formula
$$
F_n=\sum_{\sigma\in S_n}\biggl(\prod_{k=1}^n
((-1)^{k+1}s_k)^{m_k(\sigma)}\biggr),
$$
where $m_k(\sigma)$ is the number of cycles of length
$k\ge1$ in the expansion of the permutation $\sigma$ in
the product of disjoint cycles, immediately follows from
the definition of the polynomial $F_n(s_1,\ldots ,s_n)$.
Further, the number of elements of the group $S_n$ with
$m_k$ cycles of length~$k$ is equal to
$$
n!\bigg/\prod_{k=1}^n(k^{m_k}m_k!).
$$
Hence,
$$
F_n=\sum_{{\bf m}\in\Pi(n)}
\prod_{k=1}^n\biggl((-1)^{k+1}\frac{s_k}k\biggr)^{m_k},
$$
where $\Pi(n)$ denotes the set of partitions of the
number $n$ and ${\bf m}=(m_1,\ldots ,m_n)$ is a partition
of the number~$n$ such that $m_k$~is the multiplicity of
the number~$k$ in this partition. Thus,
$\sum_{k=1}^nkm_k=n$. Hence,
$$
F_nt^n=\sum_{{\bf m}\in\Pi(n)}
\prod_{k=1}^n\biggl((-1)^{k+1}\frac{s_kt^k}k\biggr)^{m_k},
$$
and we obtain the following assertion.
\medskip

{\bf Theorem~13}
$$
\sum_{n=0}^\infty F_n\frac{t^n}{n!}
=\exp\biggl(\sum_{k=1}^\infty(-1)^{k+1}s_k\frac{t^k}k\biggr).
$$

This result implies the following recursion.

\medskip

{\bf Proposition 14}
$$
F_n=(n-1)!\sum_{k=1}^n(-1)^{k+1}s_k\frac{F_{n-k}}{(n-k)!}\,.
$$

{\bf Proof} We set
$$
F(t)=\sum_{n=0}^\infty F_n\frac{t^n}{n!}\,.
$$
Then Theorem~13 implies the relation
$$
F'(t)=F(t)(s_1-s_2t+s_3t^2-\ldots),
$$
which immediately implies the desired formula. 

Let us now give a characterization of the polynomials
$F_n$ using differential operators.

\medskip

{\bf Lemma 15} 
{\it Let $\displaystyle{d=\sum_{r=2}^\infty
rs_{r-1}\frac\partial{\partial s_r}}$. Then the following
assertions hold\/:
\begin{itemize} 
\item[(a)] 
$\displaystyle{\frac{\partial F_n}{\partial s_1}=nF_{n-1}}$;

\item[(b)] 
$dF_n=-n(n-1)F_{n-1}$;

\item[(c)] 
$\displaystyle{\left[\frac{\partial}{\partial s_k},d \right]
=(k+1)\frac{\partial}{\partial s_{k+1}}}$;

\item[(d)] 
$\displaystyle{(\operatorname{Ker}\frac{\partial}{\partial
s_1}) \cap(\operatorname{Ker}d)}$ consists entirely of 
constants. 
 \end{itemize}}

{\bf Proof} (a) Differentiating the right-hand side of
the relation in Theorem~13 with respect to $s_1$, we
obtain $\displaystyle{\frac{\partial F(t)}{\partial s_1}=tF(t)}.$
Hence, 
$$
\sum_{n=0}^\infty\frac{\partial F_n}{\partial
s_1}\frac{t^n}{n!} =t\sum_{n=0}^\infty
F_n\frac{t^n}{n!}\,.
$$
Equating the coefficients of like powers of $t$, we
obtain the result in~(a).

\medskip
(b) Applying the operator $d$ to the same relation, we
obtain
$$
dF(t)=F(t)[-s_1t^2+s_2t^3-s_3t^4+\ldots]=-t^2F'(t).
$$
Equating the coefficients of $t^n$, we obtain 
$$
\frac{dF_n}{n!}=-\frac{F_{n-1}}{(n-2)!}\,.
$$
(c) 
$\displaystyle{ d\frac\partial{\partial
s_k}-\frac\partial{\partial s_k}d =\sum_{r=2}^\infty
rs_{r-1}\frac\partial{\partial s_r}\,
\frac\partial{\partial s_k}-\frac\partial{\partial s_k}
\sum_{r=2}^\infty rs_{r-1}\frac\partial{\partial s_k}
=(k+1)\frac\partial{\partial s_{k+1}}}$.

(d) Let $\displaystyle{\frac{\partial f}{\partial s_1}=0}$. Then
$f=f(s_2,s_3,\ldots )$. If $f$ belongs to the intersection
of kernels of the operators $\displaystyle{\frac\partial{\partial
s_1}}$ and $d$, then $\displaystyle{\frac{\partial}{\partial s_2}f=0}$ by
assertion~(c). Assertion~(c) enables us to
complete the proof of (d) by induction.

Lemma~15 immediately implies the following result.

\medskip

{\bf Theorem~16} 
{\it The sequence of polynomials $\{F_n(s_1,\ldots ,s_n)\}$ is
completely characterised by the following
properties\/: 
 \begin{itemize}
\item[0)] 
$F_n(0)=0$, $n=1,2,\ldots $,

\item[1)] 
$F_0=1$,

\item[2)] 
$\displaystyle{\frac{\partial F_n}{\partial s_1}=nF_{n-1}}$,

\item[3)] 
$dF_n=-n(n-1)F_{n-1}$,
\end{itemize}
that is, the generating function
$F(t)=F(t;s_1,s_2,\ldots )$ is the unique solution of the
system of equations 
$$
\frac{\partial}{\partial s_1} F(t)=tF(t),
\quad
dF(t)=-t^2\frac{\partial}{\partial t}F(t)
$$
under the initial condition $F(0;s_1,s_2,\ldots )=1$.}

\bigskip
\pagebreak

\S 4. {\large \bf Frobenius algebras}

\bigskip
Let $A$ be an associative algebra over~$\mathbb{C}$. A
linear map $f : A\to\mathbb{C}$ is said to be {\it tracial\/}
(or {\it trace-like\/}) if $f(ab)=f(ba)$ for any
$a$ and~$b$ in~$A$.

\medskip

{\bf Definition 5} 
{\it By a {\it Frobenius algebra\/}  
we mean an algebra $A$ together with a tracial linear
map $f : A\to\mathbb{C}$ such that the bilinear form 
$$
A\times A\to\mathbb{C}, \qquad (a,b)\mapsto f(ab)
$$
is non-degenerate. }

We denote by $J(A)$ the Jordan algebra whose additive
structure is that of~$A$ and with the product
$$
a\circ b=\frac12(ab+ba).
$$
Let us consider a basis $\{e_i : i\in J\}$ for $A$ and
denote the corresponding structure constants by $a_{ij}^k$
(that is, $e_ie_j=\sum a_{ij}^ke_k$). Then the numbers
$a_{(ij)}^k=\frac12(a_{ij}^k+a_{ji}^k)$ are the
structure constants of the Jordan algebra $J(A)$ with
respect to the same basis.

\medskip

{\bf Theorem~17} 
{\it Let $(A,f)$ be a Frobenius
algebra. Then the structure constants of the Jordan
algebra $J(A)$ are completely defined by the
homomorphisms $\Phi_k=\Phi_k(f)$, $k=1$, $2$, and~$3$.}

A result related to Theorem~17 was obtained in~\cite{[12]}
under other assumptions including  an
assertion similar to the conclusion of  
Theorem~2.8 in~\cite{[3]}. Our theorem therefore seems to 
give a stronger result.

\medskip

{\bf Corollary 18} 
{\it For a given homomorphism $f : A\to\mathbb{C}$ the linear maps
$\Phi_k=\Phi_k(f) : A^{\otimes k}\to\mathbb{C}$, $k\ge4$, are
determined by the maps $\Phi_1,\Phi_2,\Phi_3$.}

Corollary~18 implies the well-known result of~\cite{[13]}
that a finite group is completely determined by
its Frobenius $k$-characters for $k=1,2,3$.

\medskip

{\bf Proof of Theorem~17} 
By definition, 
$$ \Phi_2(e_i,e_j)=f(e_i)f(e_j)-f(e_ie_j). $$
We set $R_{ij}=f(e_ie_j)$. Then 
$$ R_{ij}=\sum_ra_{ij}^rf(e_r)=f(e_i)f(e_j)-\Phi_2(e_i,e_j).
 \hspace{1in} (4.1) $$

\pagebreak

Similarly, 
$$ \Phi_3(e_i,e_j,e_k)=f(e_i)f(e_j)f(e_k)-f(e_i)\sum_ra_{jk}^rf(e_r)
$$
$$-f(e_j)\sum_ra_{ik}^rf(e_r)-f(e_k)\sum_ra_{ij}^rf(e_r)
+\sum_{r,s}(a_{ij}^r+a_{ji}^r)a_{rk}^sf(e_s). $$
Hence, 
$$ \sum(a_{ij}^r+a_{ji}^r)R_{rk}\hspace{2in}$$
$$=f(e_i)R_{jk}+f(e_j)R_{ik}+f(e_k)R_{ij}
-f(e_i)f(e_j)f(e_k)+\Phi_3(e_i,e_j,e_k)
 (4.2) $$

It follows from the definition of the Frobenius algebra
that the matrix $R_{ij}=f(e_ie_j)=\sum_ra_{ij}^rf(e_r)$
is symmetric and non-degenerate. Regarding the
equation~\thetag{4.2} for fixed $i$ and~$j$ as a system
of linear equations with respect to the vector
$(a_{(ij)}^r$, $r=1,\ldots ,n)$, we see that this system
has a unique solution. This completes the proof of the
theorem. 

A more explicit answer can be obtained for a commutative
Frobenius algebra $A$.
We set 
$$
R_i=f(e_i), \quad R_{ij}=f(e_ie_j), \quad
R_{ijk}=f(e_ie_je_k). 
$$
In terms of the values of the homomorphisms $\Phi_2$ and
$\Phi_3$ we obtain
$$
R_{ij}=\Phi_2(e_i,e_j)-R_iR_j,
\\
2R_{ijk}=\Phi_3(e_i,e_j,e_k)+R_iR_{jk}+R_jR_{ik}+R_kR_{ij}-R_iR_jR_k.
$$
The matrix $R_{ij}$ is invertible by the definition of
Frobenius algebra. We denote by $R^{ij}$ the matrix
inverse to~$R_{ij}$. Direct calculations give the
following explicit formula for the structure constants
of the algebra~$A$.

\medskip

{\bf Proposition 19} 
{\it The structure constants of a commutative Frobenius
algebra are defined by the following formula\/: 
$$
a_{ij}^k=\sum_m R_{ijm}R^{mk}.
$$}

{\bf Proof} 
Using the formula 
$$
R_{ijk}=f(e_ie_je_k)=\sum_n f(a_{ij}^ne_ne_k),
$$
we obtain 
$$
\sum_m R_{ijm}R^{mk}=\sum_{m,n} a_{ij}^nR_{nm}R^{mk}=a_{ij}^k.
$$

{\bf Proof of Corollary~18} 
It suffices to prove that,
to evaluate the values of the homomorphism $\Phi_k$, one
needs the structure constants $a_{(ij)}^k$ rather
than~$a_{ij}^k$. We use induction.

The only summand in the expansion of
$\Phi_m(e_1,\ldots ,e_m)$ (in terms of the values of the
homomorphism $f$) which cannot be immediately expressed
by using the values $\Phi_r$ with $r<m$ is 
$$
\sum_{\sigma\in S_{m-1}}f(e_1
e_{\sigma(2)}e_{\sigma(3)}\ldots e_{\sigma(m)})
=\frac1m\sum_{\sigma\in S_m}f(e_{\sigma(1)}
e_{\sigma(2)}\ldots e_{\sigma(m)}). 
$$
The multiplication in the Jordan algebra $J(A)$ may
not be associative; however, the associator can be
expressed in terms of commutators, namely, 

$$
4((a\circ b)\circ c-a\circ(b\circ c))
=abc+bac+cab+cba-abc-acb-bca-cba$$
$$=[b,ac]-[b,ca]=[b,[a,c]].
$$

Since the homomorphism $f$ is tracial, it follows that
$f([a,b])=0$, and we see that 
$$
f((a\circ b)\circ
c)=f(a\circ(b\circ c)). 
$$
Hence, on iterated products of elements of the
algebra $J(A)$ the homomorphism $f$ behaves as
if this algebra is associative. Thus, 
$$
2^{m-1}\sum_{\sigma\in S_m}f(e_{\sigma(1)}\ldots
e_{\sigma(m)}) =\sum_{\sigma \in
S_m}f(e_{\sigma(1)}\circ \ldots  \circ e_{\sigma(m)}).
$$

In the case of finite groups we obtain the following
result.

Let $G$ be a finite group with the set of elements given
by $\{g_1,\ldots ,g_n\}$. We consider the Frobenius
algebra $\mathbb{C}G$ with the structure map
$f=\frac1n\chi : \mathbb{C}G\to\mathbb{C}$, where $\chi$~is the
character of the regular representation of the group~$G$.

\medskip

{\bf Corollary 20} 
{\it The linear maps $\Phi_1=f$,
$\Phi_2$, and $\Phi_3$ determine the following parts of 
the structure
of the group algebra~$\mathbb{C}G$: 
 \begin{itemize}
\item[1)] 
$\Phi_1$ determines the identity element $e$ of
$G$;

\item[2)] $\Phi_2$ determines 
inverse elements, that is, it distinguishes the pairs
$(g_i,g_j)$ such that $g_ig_j=e$;

\item[3)] $\Phi_1$, $\Phi_2$, and $\Phi_3$ define the
structure constants of the Jordan algebra $J(\mathbb{C}G)$.
\end{itemize}}

{\bf Proof} The trace $\chi$ of the regular
representation takes the value $n$ at the identity
element~$e$ of~$G$ and the value~$0$ at any other
element, which obviously implies the assertion~1)
because $\Phi_1=f=\frac{1}{n}\chi$. Further, it follows
from (4.1) that 
$$
\Phi_2(h,g)=\begin{cases} -1&\text{if}\ hg=e,\ h\ne e;
\\
0&\text{otherwise}.
\end{cases}
$$
We thus find the  pairs of mutually inverse elements. Let
us rewrite the result thus obtained in terms of
structure constants. We order the elements of~$G$ in
such a way that $g_1=e$ and write
$g_ig_j=\sum_{k=1}^na_{ij}^k g_k$. In the case under
consideration, the quantities $a_{ij}^k$ take the values
$0$ and~$1$, and $\sum_{k=1}^{n}a_{ij}^k=1$ for any
pair~$(i,j)$.  Assertion~1) means that 
$$
a_{i1}^k=a_{1i}^k=\delta_{i k}.
$$
Assertion 2) recovers the values
$$
a_{ij}^1=-\Phi_2(g_i,g_j) \quad\text{if}\quad i+j>2.
$$
We now find the information concerning the
structure constants $a_{ij}^k$ that is determined by the
map~$\Phi_3$. By~\thetag{4.2} we have
$$
\Phi_3(g_i,g_j,g_k)=\sum_r(a_{ij}^r+a_{ji}^r)a_{rk}^1
\quad\text{if}\quad i+j+k>3.
$$
We know that the  matrix $a^1_{r,k}$ is  invertible,
which enables one to complete the proof.

\medskip

{\bf Theorem~21} 
{\it Let $G$ be a finite group and let $\chi$ be the
character of the regular representation of~$G$. The
linear maps $\chi$, $\Phi_2(\chi)$, and~$\Phi_3(\chi)$
determine the group~$G$ uniquely up to isomorphism. }

{\bf Proof} 
Using Corollary~20, we immediately see that the
homomorphisms $\chi$, $\Phi_2(\chi)$, and $\Phi_3(\chi)$
define the multiplication on the set of elements of the
group $G$ up to permutation of the factors.

The following lemma and the fact that that a group is 
isomorphic to its opposite (using the map $g\mapsto g^{-1}$) 
enables one to complete the proof of
this theorem. 

\medskip

{\bf Lemma~22 \rm (Mansfield \cite{[14]})} 
{\it Let $G$ be a finite group in which the multiplication
rule is not known precisely but the set $\{gh,hg\}$ is
known for each pair of elements $g,h\in G$. Then one can
recover the set of functions $\{m,m^{\text{{op}}} : 
G\times G\to G\}$, where $m(x,y)=xy$ and
$m^{\text{{op}}}(x,y)=yx$. }

The proof of this lemma presented in~\cite{[14]} uses an
elementary (but tricky) examination of 
cases. For completeness of the exposition, we present a
proof in Appendix~A (\S~5).

In \cite{[14]} this lemma was applied to prove the result
of E.~Formanek and D.~Sibley~\cite{[8]} that the
group determinant of a finite group (which was first
considered by Dedekind) defines a group.

The group determinant of a group $G$ is defined as
follows.

Choose a one-to-one correspondence between the
set of elements $\{g_1,\ldots , g_n\}$ of
the group~$G$, where $g_1=e,$ and the set of commuting variables
$\{x_1,\ldots ,x_n\}$.
The group determinant is the determinant of the matrix
obtained from the multiplication table of the group by
the substitution $g_i\to x_i$. The group determinant is
a homogeneous polynomial of degree~$n$ in
$\{x_1,\ldots ,x_n\}$ with integral coefficients.

Introduce the matrix $M_G=(m_{ij}=x_k)$, where
the correspondence $(i,j)\to k$ is given by
$g_ig_j^{-1}=g_k$. It is clear that $M_G$ differs from
the multiplication table only by the order of the columns, and
hence $D_G=\pm\det M_G$. We note that $\det
M_G=x_1^n+\ldots$.

Let us consider the $n$-dimensional linear space
$\mathbb{C}^n\simeq\mathbb{C}G$ and the right regular representation
of the group~$G$ in the basis $\{g_1,\ldots ,g_n\}$:
$$
T : \mathbb{C}^n\to\mathbb{C}^n, \qquad T(g)g_k=g_kg^{-1}.
$$
We denote by $T_i$ the matrix of the action of the
operator $T(g_i)$ in this basis. Then 
$$
M_G=\sum x_iT_i.
$$

Consider the character of the representation
$T$, 
$$
\chi : \mathbb{C}G\to\mathbb{C}
$$
and extend it to a linear $\mathbb{C}[x_1,\ldots ,x_n]$-homomorphism 
$$
f : \mathbb{C}G[x_1,\ldots ,x_n]\to\mathbb{C}[x_1,\ldots ,x_n].
$$
\medskip

{\bf Lemma 23} 
{\it We set $a=\sum_{i=1}^nx_ig_i\in\mathbb{C}G[x_1,\ldots ,x_n]$.
Then 
$$
\Phi_n(f)(a,\ldots ,a)=D_G.
$$}
The expansion of the group determinant into
irreducible factors obtained by Frobenius in~\cite{[10]}
can  be obtained immediately from the expansion of the
regular representation~$T$ into irreducible
representations: {\it let $\chi=\chi_1+\ldots+\chi_q$ be
an expansion of the regular representation $T$ into
irreducible representations and let $n_i$~be the
dimension of the representation with the
character~$\chi_i$; then }
$$
D_G=\prod_{i=1}^q\Phi_{n_i}(\chi_i)^{n_i}.
$$

\bigskip

\S 5.  {\large \bf  Appendix A. Proof of Mansfield's lemma}
\bigskip

As we promised, we give here a  proof of the
following result.

\medskip

{\bf Lemma 24} 
{\it Let $G$ be a group with the multiplication $(x,y)\to
xy$. In this case, if $*$~is an associative
multiplication on~$G$ related to the original
multiplication in such a way that for every pair $x,y\in
G$ the product $x*y$ is equal to either $xy$ or $yx$,
then $x*y=xy$ for any $x,y\in G$ or $x*y=yx$ for any
$x,y\in G$.}

We reformulate the lemma as follows. Let
$A=\{(x,y):x*y=xy\}$ and $B=\{(x,y):x*y=yx\}$. Then it
follows from the conditions of the lemma that $A\cup
B=G\times G$, and one must prove that $A=G\times G$
or $B=G\times G$. It is clear that the set $A\cap B$ is
symmetric.

The proof of the lemma follows Exercise 26 in~\S\,4
in~\cite{[1]}. For a given pair $x,y\in G$ we have the
following assertion.

\medskip

{\bf Fact~25} $x*y=y*x\Longleftrightarrow xy=yx$.

The implication $xy=yx\Longrightarrow x*y=y*x$ is
obvious. In the converse direction, let $x*y=y*x$. In
this case, without loss of generality one can set
$x*y=xy$.  If this is not the case, then
$y*x=x*y=yx$, and we can transpose the variables
$x$ and~$y$.

In what follows we assume that $xy\ne yx$, because
otherwise there is nothing to prove. Let us consider the
element $x*x*y=x^2*y$ and assume that it differs
from~$x^2y$. At the same time it follows from the above
that $x*x*y=x*xy$ is equal to $xyx$ or
$x^2y$. Thus, $yx^2=xyx$, and hence $yx=xy$. The
contradiction thus obtained shows that $x*x*y=x^2y$.

We have $x*y*x*y=xy*xy=(xy)^2$. On the other hand, by
assumption, the element $x*y*x*y=x*x*y*y=x^2y*y$ is
equal to $x^2y^2$ or $yx^2y$, that is, $(xy)^2=x^2y^2$
or $(xy)^2=yx^2y$. In both the cases we arrive at the
relation $xy=yx$, which is a contradiction.

\medskip

{\bf Corollary 26} 
{\it $x*y=xy\Longleftrightarrow
y*x=yx$, that is, the set $A$ is symmetric.}

{\bf Proof} Let $x*y=xy$ and $y*x=xy$. Then it follows
from Fact~25 that $xy=yx$. Thus, $x*y=xy\Longrightarrow
y*x=yx$. The converse implication follows by symmetry.

Since both the sets $A$ and $A\cap B$ are symmetric and
$A\cup B=G\times G$, we immediately see that $B$ is a
symmetric set. Thus, we have proved the following
assertion.

\medskip

{\bf Corollary 27} 

$x*y=yx\Longleftrightarrow
y*x=xy$. 

\medskip

{\bf Fact 28} {\it $x*y*x=xyx$ for every $x,y\in G$.}

{\bf Proof} Let $xy=yx$, then $x*y*x=xy*x=xyx$ or
$x^2y$. The desired relation holds because $xyx=x^2y$.
On the other hand, suppose $xy\ne yx$ and $x*y*x\ne xyx$.
Then only the following two cases are possible.

a) $x*y=xy$ and $y*x=yx$. Then $x*y*x=xy*x=x^2y$ or
$xyx$. Since the multiplication is associative, it
follows that $x*y*x=x*yx$ is equal to $xyx$ or $yx^2$.
Thus, $x*y*x=x^2y=yx^2$. Moreover,
$x*x*y=x^2*y=x^2y=x*y*x$. Hence, $x*y=y*x$, which
contradicts our assumptions.

b) $x*y=yx$ and $y*x=xy$. We again have
$x*y*x=yx*x=yx^2$ and $x* y*
x=x* xy=x^2y$. Thus, $x*x*y=x*y*x$ and
$x*y=y*x$. This contradiction completes the
proof of Fact~28. 

We consider $G\times G$ as a square array with rows
$R_x=\{(x,y):y\in G\}$.

\medskip

{\bf Fact~29} 
{\it For each $x\in G$ we have either
$R_x\subset A$ or $R_x\subset B$. }

{\bf Proof} 
Suppose that Fact~29 fails. Then there is an element $x$
such that $R_x$ intersects both the sets $G\times
G\setminus A$ and $G\times G\setminus B$, that is, there
are elements $y,z\in G$ such that $x*y=yx\ne xy=y*x$ and
$x*z=xz\ne zx=z*x$. Let us show that under these
assumptions we have $yxz=zxy$. Two cases are possible
here.

a) $y*z\ne z*y$. As was shown above, in this case $yz\ne
zy$. The element $z*x*y$ is equal to either $zxy$ or
$yzx$. This element is also equal to $z*yx=zyx$ or
$yxz$. However, $yzx\ne zyx$, $yzx\ne yxz$ and $zxy\ne
zyx$, and thus $z*x*y=zxy=yxz$.

b) $y*z=z*y$ $(=yz=zy)$. We show that the assumption
$ \; yxz\ne zxy$ leads to a contradiction. We consider
the element $z*x*y$. It is equal to $zx*y$ which equals $yzx$ or 
$zxy$, but it is  also  equal to $z*yx=zyx$ or $yxz$. Under
our assumptions and using that $zxy\ne zyx$ and $yzx\ne yxz$, 
we obtain $z*x*y=yzx=zyx$. Now we consider the element
$x*z*y$. It differs from $z*x*y=yzx$. However, it is
equal to $x*y*z=x*yz=xyz$ or to $yzx$. Hence,
$x*z*y=xyz=xzy$. Finally, we consider the element
$z*y*x$. It differs from $z*x*y=yzx$ and is equal to
$yz*x=xyz$ or  $yzx$. Hence, $z*y*x=xyz=xzy$. Thus,
$x*(z*y)=(z*y)*x$, that is, $x*zy=zy*x$, and therefore
$xzy=zyx$. This implies that $z*x*y=x*z*y$, and thus
contradicts the condition $x*z\ne z*x$.

Let us finally consider the element $x*z*x*y$. According to
Fact~28, it is equal to $xzx*y=xzxy$ or $yxzx$. On the
other hand, it is equal to $xz*yx=xzyx$ or $yxxz$. We
must again consider two cases.

a) $x*z*x*y=xzxy$ is equal to $xzyx$ or $yxxz$. In the
first case we obtain $xy=yx$.
The other case, together with the above calculations,
gives $yxxz=xzxy=xyxz$, which gives $xy=yx$. In both cases 
we get that $xy=yx$ which is a contradiction.

b) $x*z*x*y=yxzx$ is equal to $xzyx$ or $yxxz$. In the
second case we have $xz=zx$, that is, a contradiction.
Thus, $xzyx=yxzx=zxyx$, which also gives $xz=zx$. This
completes the proof of Fact~29.

According to Fact~29, for every $x$ we have either $R_x\subset A$
or $R_x\subset B$. Suppose that $R_x\subset A$ but
$R_x\not\subset B$, and that $R_u\subset B$ but
$R_u\not\subset A$. Then $x*y=xy$ for any $y$, and there
is a~$z$ such that $xz\ne zx$. Similarly, $u*v=vu$ for
any~$v$, and there is a~$w$ such that $uw\ne wu$.

Consider the element $x*u=xu$. It is equal to
$ux=u*x$, and $R_{xu}$ is a subset of $A$ or $B$.
Suppose that $R_{xu}\subset A$. We then consider the
element $x*u*w=xu*w=xuw$, which must also be equal to
$x*wu=xwu$, which is a contradiction. We now assume that
$R_{xu}\subset B$ and consider the element
$z*x*u=z*xu=xuz=uxz$, which is also equal to $zx*u=uzx$.
This implies that $xz=zx$, a contradiction. This proves
Mansfield's lemma.

\bigskip

\S 6.  {\large \bf  Appendix B. The algebra of multi-symmetric
polynomials}

\bigskip
Let $n>1$. Recall that a polynomial $p(x_1,\ldots ,x_n)$, where

\noindent
$x_k=(x_{1k},\ldots ,x_{mk})\in\mathbb{C}^m$, is said to be
multi-symmetric if it is invariant with respect to all permutations of
the set of arguments $(x_1,\ldots ,x_n)$. The algebra of polynomial
functions on the algebraic variety
$\operatorname{Sym}^n(\mathbb{C}^m)$ can be canonically identified
with the algebra of multi-symmetric polynomials, which we denote by
$\mathcal{SP}^n(\mathbb{C}^m)$. For any set of non-negative integers
$\omega=(i_1,\ldots ,i_m)$ we introduce the multi-symmetric Newton
polynomials $p_{\omega}$,
$$
p_{\omega}(x_1,\ldots ,x_n)=\sum_{k=1}^{n}x_{1k}^{i_1}\ldots
x_{mk}^{i_m};
$$
for $m=1$ these are the classical
Newton polynomials.

We set $|\omega|=i_1+\ldots +i_m$. Let $\{z_{\omega}\}$, $\omega\in\mathbb
{Z}_{\ge0}^m$, be a set of commuting variables graded by the rule $\deg
z_{\omega}=|\omega|$.  We consider the graded polynomial ring $L=\mathbb
{C}[z_\omega]$. We introduce in~$L$ a system of homogeneous polynomials
$\{\mathcal{F}_{\omega_1,\ldots ,\omega_j}\}$,
$\deg \mathcal{F}_{\omega_1,\ldots ,\omega_j}=|\omega_1|+\ldots +|\omega_j|$,
by the ``Frobenius-type'' formulas,

$$ \mathcal{F}_{\omega_1}=z_{\omega_1},\quad
\mathcal{F}_{\omega_1,\omega_2}=z_{\omega_1}z_{\omega_2}-z_{\omega_1+\omega_2}
$$  and by recurrence:  
$$\mathcal{F}_{\omega_1,\ldots ,\omega_{j+1}}
=z_{\omega_1}\mathcal{F}_{\omega_2,\ldots ,\omega_{j+1}}
-\mathcal{F}_{\omega_1+\omega_2,\omega_3,\ldots ,\omega_{j+1}} -\ldots
-\mathcal{F}_{\omega_2,\ldots ,\omega_j,\omega_1+\omega_{j+1}}. $$
Following the above scheme (see \S~3), one can readily
obtain an explicit expression for the polynomial
$\mathcal{F}_{\omega_1, \ldots ,\omega_j}$. In what
follows it is important to note that the linear part of the
polynomial $\mathcal{F}_{\omega_1, \ldots ,\omega_j}$ is
equal to 

\noindent
$(-1)^{j-1}(j-1)! z_{\omega_1+\ldots +\omega_j}$.

The fact that the ring $L$ can be canonically 
identified with the polynomial ring on the linear space
$\operatorname{Hom}(\mathbb{C}[u_1,\ldots ,u_m],\mathbb{C})$
enables one to pass to the homomorphism of polynomial
rings which is induced by the evaluation map~$\mathcal{ E}$
and we obtain the following result from Theorem 11.

\medskip

{\bf Theorem 30} 
{\it For given $n$ and $m$ the ring homomorphism
$$ \mathcal{E}^* : L \to \mathcal{SP}^n(\mathbb{C}^m) $$
is an epimorphism.}

If $\mathcal{SYZ}(m,n)$ is the kernel of the
homomorphism $\mathcal{E}^*$, then $\mathcal{SYZ}(m,n)$ is
the ideal of the ring $L$ generated by the polynomials
$\mathcal{F}_{\omega_1,\ldots ,\omega_{n+1}}$ such that
$|\omega_1|,\ldots,|\omega_{n+1}|> 0$.

Consider the sub-ring  $L_n\subset L$ generated by
the elements $\{z_\omega\}$ with $|\omega|<n+1$. Using
the above form of the linear part of the polynomial
$\mathcal{ F}_{\omega_1,\ldots ,\omega_{n+1}}$, we obtain the
following assertion.

\medskip

{\bf Corollary 31} 
{\it The restriction of the homomorphism $\mathcal{E}^{*}$ to
$L_n$ gives an epimorphism 
$$
\mathcal{ E}^*_n : L_n\to\mathcal{SP}^n(\mathbb{C}^m).
$$
We denote by $\operatorname{Syz}(m,n)$ the kernel of the
homomorphism $\mathcal{E}^*_n$. Then
$$
\operatorname{Syz}(m,n)=\mathcal{SYZ}(m,n)\cap L_n.
$$}

\noindent
By construction, $L_n$ is the polynomial ring on $\mathbb{C}^N$, where
$N= \left(\begin{array}{c}n+m\\ n \end{array}\right) - 1.$ We obtain the following assertion.

\medskip

{\bf Corollary 32} 
$$
\operatorname{Sym}^n(\mathbb{C}^m)\sim\operatorname{Spec}
\bigl(L_n/\operatorname{Syz}(n,m)\bigr)\subset\mathbb{C}^N.
$$

In conclusion we note that the assertion claiming that
the homomorphism $\mathcal{E}^*$ is onto is
equivalent to the first fundamental theorem of the
classical invariant theory claiming that a
multi-symmetric polynomial in $n$ vector arguments can
be expressed (but not uniquely) as a polynomial
in the multi-symmetric Newton polynomials $p_\omega$,
$|\omega|\le n$. The classical case $m=1$ is the only
case in which this expression is unique.

V.~A.~Steklov Mathematical Institute,

Russian Academy of Sciences;

buchstab@mendeleevo.ru

\medskip
School of Mathematics,

University of Edinburgh;

E.Rees@ed.ac.uk

\end{document}